\documentclass{article}
\usepackage{amsmath}

\usepackage{amssymb,latexsym,amsmath}
\usepackage{color}

\usepackage{latexsym,amssymb}
\usepackage{hyperref}
\usepackage{graphicx}
\usepackage[frenchb,english]{babel}
\usepackage[matrix,arrow]{xy}
\def \be{\begin{eqnarray*}}
\def \ee{\end{eqnarray*}}
\def \ben{\begin{enumerate}}
\def \een{\end{enumerate}}
\def \beit{\begin{itemize}}
\def \eeit{\end{itemize}}
\def \bui#1#2{\mathrel{\mathop{\kern 0pt#1}\limits^{#2}}}
\def \buil#1#2{\mathrel{\mathop{\kern 0pt#1}\limits_{#2}}}
\def \bfll{\begin{flushleft}}
\def \efll{\end{flushleft}}
\def \bflr{\begin{flushright}}
\def \eflr{\end{flushright}}

\newtheorem{ethm}{Theorem}[section]

\newtheorem{ecor}[ethm]{Corollary}
\newtheorem{prop}[ethm]{Proposition}

\newcommand{\R}{\mathbb{R}}

\title{Riemannian foliations with parallel or harmonic basic forms}
\author{Fida EL Chami\footnote{Lebanese University, Faculty of Sciences II, Department of Mathematics, P.O. Box 90656 Fanar-Matn, Lebanon,
E-mail: \texttt{fchami@ul.edu.lb}},\, Georges Habib\footnote{Lebanese University, Faculty of Sciences II, Department of Mathematics, P.O. Box 90656 Fanar-Matn, Lebanon,
E-mail: \texttt{ghabib@ul.edu.lb}},\, Roger Nakad \footnote{Notre Dame University-Louaiz\'e, Faculty of Natural and Applied Sciences, Department of Mathematics and Statistics, P.O. Box 72, Zouk Mikael, Lebanon, E-mail: \texttt{rnakad@ndu.edu.lb}}}

\begin{document}
\date{}
\maketitle
\begin{abstract} 
\noindent In this paper, we consider a Riemannian foliation whose normal bundle carries a parallel or   harmonic basic form. We estimate the norm of the O'Neill tensor in terms of the curvature data of the whole manifold. Some examples are then given. 
\end{abstract}\vspace{+0.5cm}
{\it $2010$ Mathematics subject classifications:}  53C12, 53C20, 53C24, 57R30.

\section{Introduction}
In \cite{Gro}, J. F. Grosjean obtained some non-existence results on minimal submanifolds carrying parallel or harmonic forms. Indeed, given a Riemannian manifold $(M^m, g)$ admitting a parallel $p$-form and let $(N^n, h)$ be a Riemannian manifold satisfying a certain curvature pinching condition depending on $m$ and $p$, he proved that there is no minimal
immersion from $M$ into $N$. His proof is based on computing the curvature term (which is zero in this case) in the Bochner Weitzenb\"ock formula and using the Gauss formula relating the curvatures of $M$ and $N$. As a consequence, he deduced various rigidity results when $N$ is the hyperbolic space $\mathbb{H}^n,$  the Riemannian product $\mathbb{H}^r \times \mathbb{S}^s$ or the complex hyperbolic space $\mathbb{CH}^n.$\\

\noindent In the same spirit, he proved that for any compact manifold $(M^m, g)$ carrying a harmonic $p$-form (or a non-zero $p$th betti number $b_p(M)$)
and isometrically immersed into a Riemannian manifold $(N^n,h)$, there exists at least a point 
$x$ of $M$ so that (see also \cite{S})
$$\frac{m}{\sqrt{p}} \left(\frac{p-1}{p}\right) \vert B(x)\vert \vert H(x) \vert \geq k(x)- \left(\frac {p-1}{p}\right) ((m-1)\bar K_1 + \bar \rho_1) (x),$$
and 
$$m \left(\frac{p-1}{\sqrt{p}} + \frac{m-p-1}{\sqrt{m-p}}\right) \vert B(x)\vert \vert H(x)\vert \geq \mathrm{Scal}^M(x) - (m-2)[(m-1)\bar K_1 + \bar \rho_1] (x),$$
where $|B(x)|$, $ H(x)$, $k(x)$ and $\bar \rho_1(x)$ denote respectively the norm of the second fundamental
form $B,$ the mean curvature of the immersion, the smallest eigenvalue of the Ricci curvature of $M,$ the largest
eigenvalue of the curvature operator of $(N^n, h)$ and $\bar K_1
(x)$ is the largest sectional curvature of $N$. These inequalities come from a lower bound of the curvature term (which is non-positive at the point $x$) in the Bochner Weitzenb\"ock formula. Thus, if the manifold $(M^m, g)$
is minimally immersed into $(N^n, h)$ and satisfying the pinching condition
$$\mathop\mathrm{min}\limits_M (\mathrm{Scal}^M) > (m-2) \big ( (m-1) \mathop\mathrm{max}\limits_N (\bar K_1) +  \mathop\mathrm{max}\limits_N (\bar \rho_1) \big), $$
then $(M^m, g)$ is a sphere of homology (see also \cite{L}).\\

\noindent In this paper, we investigate the study of foliated manifolds whose normal bundle carries a particular form. In fact, we consider a Riemannian manifold $(M,g)$ equipped with a Riemannian foliation $\mathcal{F},$ which roughly speaking, is the decomposition of $M$ into submanifolds (called leaves) given by local Riemannian submersions to a base manifold. We assume that the normal bundle of the foliation carries a parallel $p$-form (resp. harmonic), with respect to the connection defined in Section \ref{sec:1}. This corresponds locally to the existence of such a form on the base manifold of the submersions. When shifting the study from immersions to submersions, many objects are replaced by their dual. In particular, the O'Neill tensor \cite{ON} plays the role of the second fundamental form and thus, we aim to estimate  the norm of the O'Neill tensor in terms of different curvatures data of the manifold $M$. The main tool is to use the transverse Bochner Weitzenb\"ock formula for foliations \cite{HR}. Recall that this tensor completely determines the geometry of the foliation. Indeed, it vanishes if and only if the normal bundle of the foliation is integrable.\\  

\noindent The paper is organized as follows. In Section \ref{sec:1}, we recall some well-known facts on differential forms and review some preliminaries on Riemannian foliations. In Section \ref{sec:3}, we treat the case where the normal bundle carries a parallel form. We compute the curvature term in the transverse Bochner Weitzenb\"ock formula and relate it to the curvature of the manifold $M$ using the O'Neill formulas. We then deduce a lower bound estimate for the O'Neill tensor (see Thm. \ref{thm:1} for $p>1$ and Cor. \ref{cor:1} for a rigidity result when $p=1$). In the last section, we study the case where there exist a harmonic form. As before, we deduce a new estimate of the O'Neill tensor (see Thm. \ref{thm:harmo}).  

\section{Preliminaries}\label{sec:1}
\setcounter{equation}{0} Let $(M,g)$ be a Riemannian manifold of
dimension $n$ and $\nabla^M$ be the Levi-Civita connection associated with the metric $g$. In all
the paper, we make the following notations for the curvatures,
$R^M(X,Y)=\nabla^M_{[X,Y]}-[\nabla^M_X,\nabla^M_Y]$ and
$R^M_{XYZW}=g(R^M(X,Y)Z,W)$ for any $X,Y,Z, W\in \Gamma(TM)$. We will
denote respectively by $K^M_0(x)$ and $K^M_1(x)$ the smallest and
the largest sectional curvature and by $\rho^M_0(x)$ and
$\rho^M_1(x)$ the smallest and largest eigenvalue of the curvature
operator $\rho^M(X\wedge Y, Z\wedge W)=g(R^M(X,Y)Z,W)$ at a point
$x\in M.$ Thus, we have the following inequalities
\begin{equation}\label{eq:12}
\rho^M_0(x)\leq K^M_0(x)\leq K^M_1(x)\leq \rho^M_1(x).
\end{equation}
Now, let us recall some definitions on forms. The inner product of any two $p$-forms $\alpha$ and $\beta$ is defined as
$$\langle\alpha,\beta\rangle=\frac{1}{p!}\sum_{1\leq i_1 \leq \cdots \leq i_p\leq n}\alpha(e_{i_1},e_{i_2},\cdots,e_{i_p})\beta(e_{i_1},e_{i_2},\cdots,e_{i_p}),$$
where $\{e_1,\cdots,e_n\}$ is an orthonormal frame of $TM.$ The interior product of a $p$-form $\alpha$ with a vector field $X$ is a $(p-1)$-form defined by
$$(X\lrcorner\alpha)(X_1,\cdots,X_{p-1})=\alpha(X,X_1,\cdots,X_{p-1}).$$ More generally, the interior product of $\alpha$ with $s$ vector fields $X_1, X_2,\cdots,X_s$ is a $(p-s)$-form which is defined as the following
$$((X_1\wedge \cdots\wedge X_s)\lrcorner\alpha)(Y_1,\cdots,Y_{p-s})=\alpha(X_s,\cdots,X_1,Y_{1},\cdots,Y_{p-s}).$$ As a consequence from the definition, we get the rule
$X\lrcorner(\omega\wedge\theta)=(X\lrcorner\omega)\wedge\theta+(-1)^p\omega\wedge(X\lrcorner \theta),$ where $p$ is the degree of $\omega.$ If the manifold is orientable, the Hodge operator $*$ defined on a $p$-form $\alpha$ satisfies the following property: 
\begin{equation}\label{eq:112}
X\lrcorner (*\alpha)=(-1)^p*(X^*\wedge \alpha).
\end{equation} 
Assume now that $(M^n,g)$ is endowed with a Riemannian foliation $\mathcal{F}$ of codimension $q$. That means $\mathcal{F}$ is given by an integrable subbundle $L$ of $TM$ of rank $n-q$ such that the metric $g$ satisfies the holonomy-invariance condition on the normal vector bundle $Q=TM/L;$ that is
$\mathcal{L}_Xg|_Q=0$ for all $X\in \Gamma(L),$ where $\mathcal{L}$ denotes the Lie derivative \cite{Rein}. We call $g$ a bundle-like metric. This latter condition gives rise to a transverse Levi-Civita connection on $Q$ defined by \cite{T}
$$
\nabla _{X} Y =
\left\{\begin{array}{ll}
\pi [X,Y],\,\,\textrm {if $X\in\Gamma(L)$},\\\\
\pi (\nabla_{X}^{M}Y),\,\,\textrm{if  $X \in \Gamma (Q)$},
\end{array}\right.
$$
where $\pi:TM\rightarrow Q$ is the projection. A fundamental property of the connection $\nabla$ is that it is flat along the leaves, that is $X\lrcorner R^\nabla=0$ for any $X\in \Gamma(L).$ Thus, we can associate to $\nabla$ all the curvatures data such as the transverse Ricci curvature ${\rm Ric}^\nabla$ and transverse scalar curvature ${\rm Scal}^\nabla.$
A basic form $\alpha$ on a Riemannian foliation is a differential form which depends locally on the transverse variables; that is satisfying the rules $X\lrcorner \alpha=0$ and $X\lrcorner d\alpha=0$ for any $X\in \Gamma(L).$ It is easy to see that the exterior derivative $d$ preserves the set of basic forms and its restriction to this set will be denoted by $d_b.$ We let $\delta_b$ the formal adjoint of $d_b$ with respect to the $L^2$-product. Then we have
$$d_b=\sum_{i=1}^q e_i\wedge\nabla_{e_i},\,\, \delta_b=-\sum_{i=1}^q e_i\lrcorner \nabla_{e_i}+\kappa\lrcorner,$$
where $\{e_i\}_{i=1,\cdots,q}$ is an orthonormal frame of $\Gamma(Q)$ and $\kappa=\mathop\sum\limits_{s=1}^{n-q}\pi(\nabla^M_{V_s}V_s)$ is the mean curvature field of the foliation, which is assumed to be a basic $1$-form. Here $\{V_s\}_{s=1,\cdots,n-q}$ is an orthonormal frame of $\Gamma(L).$ The basic Laplacian is defined as $\Delta_b=d_b\delta_b+\delta_b d_b.$ Recall that when the foliation is transversally orientable, the basic Hodge operator $*_b$ is defined on the set of basic $p$-forms as being 
$$*_b \alpha=(-1)^{(n-q)(q-p)}*(\alpha\wedge \chi_\mathcal{F}),$$ 
where $\chi_\mathcal{F}$ is the volume form of the leaves. It is also a basic $(q-p)$-form which satisfies the same property as \eqref{eq:112}. In \cite{HR}, the authors define a new twisted exterior derivative $\tilde{d}_b:=d_b-\frac{1}{2}\kappa\wedge$ and prove that the associated twisted Laplacian $\tilde{\Delta}_b:=\tilde{d}_b\tilde{\delta}_b+\tilde{\delta}_b\tilde{d}_b$ commutes with the basic Hodge operator. In particular, this shows that the twisted cohomology group (i.e. the one associated with $\tilde{d}_b$) satisfies the Poincar\'e duality. Here $\tilde{\delta}_b:=\delta_b-\frac{1}{2}\kappa\lrcorner$ denotes the $L^2$-adjoint of $\tilde{d}_b$. Moreover, they state the transverse Bochner-Weitzenb\"ock formula for $\tilde{\Delta}_b$
$$\tilde{\Delta}_b\alpha=\nabla^*\nabla\alpha+\frac{1}{4}|\kappa|^2\alpha+R(\alpha),$$
where $R(\alpha)=-\mathop\sum\limits_{j=1}^q e_j\wedge (e_i\lrcorner R^{\nabla}(e_i,e_j)\alpha).$ As for ordinary manifolds \cite{Gro}, the scalar product of $R(\alpha)$ by $\alpha$ gives after the use of the first Bianchi identity that
\begin{equation}
\langle R(\alpha),\alpha\rangle=\sum_{1\leq i,j\leq q}{\rm Ric}^\nabla_{ij}\langle e_i\lrcorner \alpha,e_j\lrcorner\alpha\rangle-\frac{1}{2}\sum_{1\leq i,j,k,l\leq q}R^\nabla_{ijkl}\langle e_j\wedge e_i\lrcorner\alpha,e_l\wedge e_k\lrcorner\alpha\rangle.
\label{eq:2}
\end{equation}
On the other hand, the geometry of a Riemannian foliation can be interpreted in terms of the so-called the O'Neill tensor \cite{ON}. It is a $2$-tensor field given for all $X,Y\in \Gamma(TM)$ by
$$A_XY=\pi^\perp(\nabla^M_{\pi(X)}\pi(Y))+\pi(\nabla^M_{\pi(X)}\pi^\perp(Y)),$$
where $\pi^\perp$ denotes the projection of $TM$ onto $L.$ By the bundle-like condition, the O'Neill tensor is a skew-symmetric tensor with respect to the vector fields $Y,Z\in \Gamma(Q)$ and it is equal to $A_Y Z=\frac{1}{2}\pi^\perp([Y,Z])$ and for any $V\in \Gamma(L)$ we have $g(A_Y V,Z)=-g(V,A_Y Z).$ Thus we deduce that the normal bundle is integrable if and only if the O'Neill tensor vanishes. If moreover the bunlde $L$ is totally geodesic, the foliation is isometric to a local product. \\\\
We point out that the curvature of $M$ can be related to the one on the normal bundle $Q$ via the O'Neill tensor by the formula \cite{ON}
\begin{equation} \label{eq:gh}
R^M_{XYZW}=R^Q_{XYZW}-2g(A_XY,A_ZW)+g(A_YZ,A_XW)+g(A_ZX,A_YW),
\end{equation}
where $X,Y,Z,W$ are vector fields in $\Gamma(Q)$. One can easily see by \eqref{eq:gh} that the norm of the O'Neill tensor $|A|^2:=\displaystyle\sum_{\substack{1\leq i\leq q\\1\leq s\leq n-q}}|A_{e_i}V_s|^2$ can be bounded at any point by 
$${\rm Scal}^\nabla-q(q-1)K^M_1\leq 3|A|^2\leq {\rm Scal}^\nabla-q(q-1)K^M_0.$$ 
In particular, if the transversal scalar curvature does not belong to the interval $[q(q-1)K^M_0,q(q-1)K^M_1]$, the normal bundle cannot be integrable.    
\section{Foliations with parallel basic forms} \label{sec:3}
\setcounter{equation}{0}
In this section, we discuss the case where the normal bundle of a Riemannian foliation carries a parallel basic form. That is a $p$-form $\alpha$ satisfying $\nabla\alpha=0.$\\

\begin{prop} \label{pro:1}Let $(M,g,\mathcal{F})$ be a Riemannian manifold with a Riemannian foliation $\mathcal{F}$ of codimension $q$. Assume that the normal bundle carries a parallel $p$-form $\alpha$. Then we have
\begin{eqnarray}\label{eq:7}
0&\leq &-\sum_{1\leq i,j,l\leq q}R^M_{lilj}\langle e_i\lrcorner\alpha, e_j\lrcorner\alpha\rangle+\frac{1}{2}\sum_{1\leq i,j,k,l\leq q} R^M_{ijkl}\langle e_j\wedge e_i\lrcorner \alpha,e_l\wedge e_k\lrcorner\alpha\rangle \nonumber\\&&+\sum_{s=1}^{n-q}\left\{|\sum_{i=1}^q A_{e_i}V_s\wedge e_i\lrcorner\alpha|^2-2\sum_{i=1}^q|A_{e_i}V_s\lrcorner\alpha|^2\right\},
\end{eqnarray}
where $\{e_i\}_{i=1,\cdots,q}$ and $\{V_s\}_{s=1,\cdots,n-q}$ are respectively orthonormal frames of $\Gamma(Q)$ and $\Gamma(L).$
\end{prop}
{\bf Proof.} From Equation \eqref{eq:gh}, we have the following formulas
\begin{eqnarray}
R^\nabla_{ijkl}=R^M_{ijkl}+2g(A_{e_i}e_j,A_{e_k}e_l)-g(A_{e_j}e_k,A_{e_i}e_l)-g(A_{e_k}e_i,A_{e_j}e_l),
\label{eq:21}
\end{eqnarray}
and that,
\begin{eqnarray}
{\rm Ric}^\nabla_{ij}=\sum_{l=1}^q\left\{R^M_{lilj}+2g(A_{e_l}e_i,A_{e_l}e_j)-g(A_{e_i}e_l,A_{e_l}e_j)-g(A_{e_l}e_l,A_{e_i}e_j)\right\}.
\label{eq:22}
\end{eqnarray}
The existence of a parallel form $\alpha$ implies that $\langle R(\alpha),\alpha\rangle=0.$ Thus plugging these last two equations into \eqref{eq:2}, we get that
\begin{eqnarray} \label{eq:3}
0&=&\sum_{1\leq i,j,l\leq q} R^M_{lilj}\langle e_i\lrcorner\alpha,e_j\lrcorner\alpha\rangle+3g(A_{e_l}e_i,A_{e_l}e_j)\langle e_i\lrcorner\alpha,e_j\lrcorner\alpha\rangle\nonumber\\
&&\sum_{1\leq i,j,k,l\leq q}\{-\frac{1}{2}R^M_{ijkl}\langle e_j\wedge e_i\lrcorner\alpha,e_l\wedge e_k\lrcorner\alpha\rangle-g(A_{e_i}e_j,A_{e_k}e_l)\langle e_j\wedge e_i\lrcorner\alpha,e_l\wedge e_k\lrcorner\alpha\rangle \nonumber\\
&&+\frac{1}{2}g(A_{e_j}e_k,A_{e_i}e_l)\langle e_j\wedge e_i\lrcorner\alpha,e_l\wedge e_k\lrcorner\alpha\rangle+\frac{1}{2}g(A_{e_k}e_i,A_{e_j}e_l)\langle e_j\wedge e_i\lrcorner\alpha,e_l\wedge e_k\lrcorner\alpha\rangle\}.\nonumber\\
\end{eqnarray}
The last two summations in the above equality are in fact equal. Indeed, using that the O'Neill tensor is antisymmetric, we find
\begin{eqnarray*}
\sum_{1\leq i,j,k,l\leq q}g(A_{e_k}e_i,A_{e_j}e_l)\langle e_j\wedge e_i\lrcorner\alpha,e_l\wedge e_k\lrcorner\alpha\rangle&=&-\sum_{1\leq i,j,k,l\leq q}g(A_{e_i}e_k,A_{e_j}e_l)\langle e_j\wedge e_i\lrcorner\alpha,e_l\wedge e_k\lrcorner\alpha\rangle\\
&=&-\sum_{1\leq i,j,k,l\leq q}g(A_{e_i}e_l,A_{e_j}e_k)\langle e_j\wedge e_i\lrcorner\alpha,e_k\wedge e_l\lrcorner\alpha\rangle\\
&=&\sum_{1\leq i,j,k,l\leq q}g(A_{e_i}e_l,A_{e_j}e_k)\langle e_j\wedge e_i\lrcorner\alpha,e_l\wedge e_k\lrcorner\alpha\rangle.
\end{eqnarray*}
On the other hand, we have
\begin{eqnarray} \label{eq:4}
\sum_{1\leq i,j,k,l\leq q}g(A_{e_i}e_j,A_{e_k}e_l)\langle e_j\wedge e_i\lrcorner\alpha,e_l\wedge e_k\lrcorner\alpha\rangle&=&\sum_{\substack{1\leq i,j,k,l\leq q\nonumber\\ 1\leq s\leq n-q}}g(A_{e_i}e_j,V_s)g(A_{e_k}e_l,V_s)\langle e_j\wedge e_i\lrcorner\alpha,e_l\wedge e_k\lrcorner\alpha\rangle\nonumber\\
&=&\sum_{\substack{1\leq i,j,k,l\leq q\\1\leq s\leq n-q}}g(A_{e_i}V_s,e_j)g(A_{e_k}V_s,e_l)\langle e_j\wedge e_i\lrcorner\alpha,e_l\wedge e_k\lrcorner\alpha\rangle\nonumber\\
&=&\sum_{s=1}^{n-q}\langle(\sum_{i=1}^q A_{e_i}V_s \wedge e_i)\lrcorner\alpha,(\sum_{k=1}^{q} A_{e_k}V_s \wedge e_k)\lrcorner\alpha\rangle\nonumber\\
&=&\sum_{s=1}^{n-q}|(\sum_{i=1}^q A_{e_i}V_s \wedge e_i)\lrcorner\alpha|^2.
\end{eqnarray}
Also we have that
\begin{eqnarray} \label{eq:5}
\sum_{1\leq i,j,l\leq q}g(A_{e_l}e_i,A_{e_l}e_j)\langle e_i\lrcorner\alpha,e_j\lrcorner\alpha\rangle &=&\sum_{\substack{1\leq i,j,l\leq q\\1\leq s\leq n-q}}g(A_{e_l}e_i,V_s)g(A_{e_l}e_j,V_s)\langle e_i\lrcorner\alpha,e_j\lrcorner\alpha\rangle\nonumber\\
&=&\sum_{\substack{1\leq i,j,l\leq q\\1\leq s\leq n-q}}g(A_{e_l}V_s,e_i)g(A_{e_l}V_s,e_j)\langle e_i\lrcorner\alpha,e_j\lrcorner\alpha\rangle\nonumber\\
&=&\sum_{\substack{1\leq l\leq q\\1\leq s\leq n-q}}\langle A_{e_l}V_s\lrcorner\alpha,A_{e_l}V_s\lrcorner\alpha\rangle=\sum_{\substack{1\leq l\leq q\\1\leq s\leq n-q}}|A_{e_l}V_s\lrcorner\alpha|^2.\nonumber\\
\end{eqnarray}
In order to estimate the last term in \eqref{eq:3}, we introduce the $p$-tensor 
\begin{eqnarray*}
\mathcal{B}^{+}(\alpha)(X_1,\cdots,X_p)=\mathop\sum\limits_{i=1}^q (e_i\lrcorner\alpha\wedge A_{e_i})(X_1,\cdots,X_p),
\end{eqnarray*}
for any $X_1,\cdots, X_p\in \Gamma(Q).$ We now proceed the computation as in \cite{Gro}. The norm of the tensor $\mathcal{B}^+(\alpha)$ is equal to
\begin{eqnarray*}
|\mathcal{B^+}(\alpha)|^2&=&\frac{1}{p!}\sum_{1\leq i_1,\cdots,i_p,i,j\leq q}\langle(e_i\lrcorner\alpha\wedge A_{e_i})_{i_1,\cdots,i_p},(e_j\lrcorner\alpha\wedge A_{e_j})_{i_1,\cdots,i_p}\rangle\\
&=&\frac{1}{p!}\sum_{\substack{1\leq i_1,\cdots,i_p,i,j\leq q\\r,t}}(-1)^{r+t}g(A_{e_i}e_{i_r},A_{e_j}e_{i_t})\alpha_{ii_1,\cdots,{\hat i_r},\cdots,i_p}\alpha_{ji_1,\cdots,{\hat i_t},\cdots,i_p}\\
&=&\frac{1}{p!}\sum_{\substack{1\leq i_1,\cdots,i_p,i,j\leq q\\r=t}}g(A_{e_i}e_{i_r},A_{e_j}e_{i_r})\alpha_{ii_1,\cdots,{\hat i_r},\cdots,i_p}\alpha_{ji_1,\cdots,{\hat i_r},\cdots,i_p}\\
&&+\frac{1}{p!}\sum_{\substack{1\leq i_1,\cdots,i_p,i,j\leq q\\r<t}}(-1)^{r+t}g(A_{e_i}e_{i_r},A_{e_j}e_{i_t})\alpha_{ii_1,\cdots,{\hat i_r},\cdots,i_p}\alpha_{ji_1,\cdots,{\hat i_t},\cdots,i_p}
\\&&+\frac{1}{p!}\sum_{\substack{1\leq i_1,\cdots,i_p,i,j\leq q\\r>t}}(-1)^{r+t}g(A_{e_i}e_{i_r},A_{e_j}e_{i_t})\alpha_{ii_1,\cdots,{\hat i_r},\cdots,i_p}\alpha_{ji_1,\cdots,{\hat i_t},\cdots,i_p}\\
&=&\frac{1}{(p-1)!}\sum_{1\leq i_1,\cdots,i_{p-1},i,j,k\leq q}g(A_{e_i}e_{k},A_{e_j}e_{k})\alpha_{ii_1,\cdots,i_{p-1}}\alpha_{ji_1,\cdots,i_{p-1}}
\\&&-\frac{2}{p!}\sum_{\substack{1\leq i_1,\cdots,i_p,i,j\leq q\\r<t}}g(A_{e_i}e_{i_r},A_{e_j}e_{i_t})\alpha_{ii_ti_1,\cdots,{\hat i_r},\cdots,{\hat i_t},\cdots,i_p}\alpha_{ji_ri_1,\cdots,{\hat i_r},\cdots,{\hat i}_t,\cdots,i_p}.
\end{eqnarray*}
Since we can choose $\frac{p(p-1)}{2}$ numbers $r,t$ with $r<t$ from the set $\{1,\cdots,p\},$ the last equality can be reduced to
\begin{eqnarray}
|\mathcal{B^+}(\alpha)|^2&=&\frac{1}{(p-1)!}\sum_{1\leq i_1,\cdots,i_{p-1},i,j,k\leq q}g(A_{e_i}e_{k},A_{e_j}e_{k})\alpha_{ii_1,\cdots,i_{p-1}}\alpha_{ji_1,\cdots,i_{p-1}}\nonumber\\
&&-\frac{1}{(p-2)!}\sum_{1\leq i_1,\cdots,i_{p-2},i,j,k,l\leq q}g(A_{e_i}e_{k},A_{e_j}e_{l})\alpha_{ili_1,\cdots,i_{p-2}}\alpha_{jki_1,\cdots,i_{p-2}}
\nonumber\\
&=&\sum_{1\leq i,j,k\leq q}g(A_{e_i}e_{k},A_{e_j}e_{k})\langle e_i\lrcorner\alpha,e_j\lrcorner\alpha\rangle-\sum_{1\leq i,j,k,l\leq q}g(A_{e_i}e_{k},A_{e_j}e_{l})\langle e_l\wedge e_i\lrcorner\alpha,e_k\wedge e_j\lrcorner\alpha\rangle\nonumber\\
&=&\sum_{1\leq i,j,k\leq q}g(A_{e_k}e_{i},A_{e_k}e_{j})\langle e_i\lrcorner\alpha,e_j\lrcorner\alpha\rangle+\sum_{1\leq i,j,k,l\leq q}g(A_{e_i}e_{l},A_{e_j}e_{k})\langle e_j\wedge e_i\lrcorner\alpha,e_l\wedge e_k\lrcorner\alpha\rangle.\nonumber\\
\label{eq:6}
\end{eqnarray}
Returning back to the Equation \eqref{eq:3} and after plugging Equations \eqref{eq:4}, \eqref{eq:5} and \eqref{eq:6}, we get the following
\begin{eqnarray*}
0&=&\sum_{1\leq i,j,l\leq q} R^M_{lilj}\langle e_i\lrcorner\alpha,e_j\lrcorner\alpha\rangle-\frac{1}{2}\sum_{1\leq i,j,k,l\leq q}R^M_{ijkl}\langle e_j\wedge e_i\lrcorner\alpha,e_l\wedge e_k\lrcorner\alpha\rangle\\
&&+|\mathcal{B^+}(\alpha)|^2+2\sum_{\substack{1\leq l\leq q\\1\leq s\leq n-q}}|A_{e_l}V_s\lrcorner\alpha|^2-\sum_{s=1}^{n-q}|(\sum_{i=1}^q A_{e_i}V_s \wedge e_i)\lrcorner\alpha|^2.
\end{eqnarray*}
Finally using the fact that $|\mathcal{B^+}(\alpha)|^2\geq 0,$ we deduce the desired inequality.
\hfill$\square$\\\\
For $p=1,$ we find by \eqref{eq:7} that the lowest sectional curvature $K^M_0$ should be non-positive. Hence we have 
\begin{ecor} \label{cor:1}
Let $(M,g,\mathcal{F})$ be a Riemannian manifold with positive sectional curvature and endowed with a Riemannian foliation $\mathcal{F}$. Then $M$ cannot carry a parallel basic $1$-form.  
\end{ecor}
\noindent In the following, we will treat the case $p\geq 2.$ For that, we aim to estimate each term in Inequality \eqref{eq:7}. As in \cite{Gro}, we define the basic $2$-form $\theta^{i_1,\cdots,i_{p-2}}=\frac{1}{2}\displaystyle\sum_{1\leq i,j\leq q}\alpha_{iji_1,\cdots,i_{p-2}}e_i\wedge e_j.$ Thus the second term of Inequality \eqref{eq:7} can be bounded from above by
\begin{eqnarray} \label{eq:8}
\frac{1}{2}\sum_{1\leq i,j,k,l\leq q}R^M_{ijkl}\langle e_j\wedge e_i\lrcorner\alpha, e_l\wedge e_k\lrcorner\alpha\rangle&=&\frac{2}{(p-2)!}\sum_{1\leq i_1,\cdots,i_{p-2}\leq q}\rho^M(\theta^{i_1,\cdots,i_{p-2}},\theta^{i_1,\cdots,i_{p-2}})\nonumber\\
&\leq &\frac{2}{(p-2)!}\rho^M_1\sum_{1\leq i_1,\cdots,i_{p-2}\leq q}|\theta^{i_1,\cdots,i_{p-2}}|^2=p(p-1)\rho^M_1|\alpha|^2.\nonumber\\
\end{eqnarray}
\noindent Using the Cauchy-Schwarz inequality and the fact that $|v\wedge w\lrcorner \alpha|\leq |v||w\lrcorner\alpha|$ for any vectors $v,w,$ the last term in \eqref{eq:7} is bounded by  
\begin{equation} 
|\sum_{i=1}^q A_{e_i}V_s\wedge e_i\lrcorner \alpha|^2\leq q \sum_{i=1}^q|A_{e_i}V_s\wedge e_i\lrcorner\alpha|^2\leq q \sum_{i=1}^q|A_{e_i}V_s\lrcorner\alpha|^2.
\label{eq:9}
\end{equation} 
Now we state our main result:

\begin{ethm} \label{thm:1} Let $(M,g,\mathcal{F})$ be a Riemannian manifold with a Riemannian foliation $\mathcal{F}$ of codimension $q\geq 4$. Assume that the normal bundle carries a parallel $p$-form $\alpha$ with $2\leq p\leq q-2$. Then we have
$$(q-2)|A|^2\geq K^M_0 q(q-1)-(p(p-1)+(q-p)(q-p-1))\rho^M_1.$$
\end{ethm}

\noindent {\bf Proof.} Plugging  the estimates in \eqref{eq:8} and \eqref{eq:9} into Inequality  \eqref{eq:7}, we get that
\begin{eqnarray}\label{eq:10}
0&\leq& -\sum_{1\leq i,j,l\leq q}R^M_{lilj}\langle e_i\lrcorner\alpha,e_j\lrcorner\alpha\rangle+p(p-1)\rho^M_1|\alpha|^2
+(q-2)\sum_{\substack{1\leq i\leq q\\ 1\leq s\leq n-q}}|A_{e_i}V_s\lrcorner\alpha|^2.\nonumber\\
\end{eqnarray}
Since $\alpha$ is a parallel $p$-form, the $(q-p)$-form $*_b\alpha$
is also parallel. Thus replacing $\alpha$ by $*_b\alpha$ in
\eqref{eq:10}, we find the inequality
\begin{eqnarray}\label{eq:11}
0&\leq &-\sum_{1\leq i,j,l\leq q} R^M_{lilj}\langle e_i\lrcorner(*_b\alpha),e_j\lrcorner(*_b\alpha)\rangle+(q-p)(q-p-1)\rho^M_1|\alpha|^2+(q-2)\sum_{\substack{1\leq i\leq q\\ 1\leq s\leq n-q}}|A_{e_i}V_s\wedge\alpha|^2.\nonumber\\
\end{eqnarray}
In the last term of \eqref{eq:11}, we use the equality \ref{eq:112} for the basic Hodge operator. Now the sum of Inequalities \eqref{eq:10} and \eqref{eq:11} gives the desired inequality after the use of
\begin{eqnarray}
\sum_{1\leq i,j,l\leq
q}R^M_{lilj}(\langle e_i\lrcorner\alpha,e_j\lrcorner\alpha\rangle+\langle e_i\lrcorner(*_b\alpha),e_j\lrcorner(*_b\alpha)\rangle)&=&
\sum_{1\leq i,j,l\leq q}R^M_{lilj}(\langle e_i\lrcorner\alpha,e_j\lrcorner\alpha\rangle+\langle e_i\wedge\alpha,e_j\wedge\alpha\rangle)\nonumber\\
&=&\sum_{1\leq i,j,l\leq q}R^M_{lilj}(\langle e_i\lrcorner\alpha,e_j\lrcorner\alpha\rangle+\langle e_j\lrcorner(e_i\wedge\alpha),\alpha\rangle)\nonumber\\
&=&\sum_{1\leq i,j,l\leq q}R^M_{lilj}(\langle e_i\lrcorner\alpha,e_j\lrcorner\alpha\rangle+\delta_{ij}|\alpha|^2-\langle e_i\lrcorner\alpha,e_j\lrcorner\alpha\rangle)\nonumber\\
&=&\sum_{1\leq i,l\leq q}R^M_{lili}|\alpha|^2,
\label{eq:cour}
\end{eqnarray}
which is greater than $K^M_0 q(q-1)|\alpha|^2.$
\hfill$\square$ \\\\

\noindent We point out that the theorem is of interest only if
$$K^M_0 q(q-1)-(p(p-1)+(q-p)(q-p-1))\rho^M_1>0,$$
which means by \eqref{eq:12} that the manifold $M$ is of positive sectional curvature.\\\\

\noindent{\bf Example:} Let us consider the round sphere $\mathbb{S}^{2m-1}$ equiped with the standard metric of constant curvature $1$. We denote by $\mathcal{F}$ the $1$-dimensional Riemannian fibers given by the action \cite{Gr}
$$e^{2i\pi t}(z_1,\cdots,z_m)=(e^{2i\pi\theta_1 t} z_1, e^{2i\pi\theta_2 t} z_2,\cdots, e^{2i\pi\theta_m t} z_m),$$ 
with $0<\theta_1\leq \theta_2\leq \cdots \leq \theta_m\leq 1.$ These foliations are Seifert fibrations (i.e. the fibers are compact) if and only if all $\theta_i's$ are rational and the Hopf fibration corresponds to the case where $\theta_1=\theta_2=\cdots=\theta_m=1.$ In the following, we will compute the O'Neill tensor of the foliation $\mathcal{F}$ and study the optimality of the estimate in Theorem \ref{thm:1}. Without loss of generality, we can assume that $\theta_1=1.$ The vector $X$ that generates $\mathcal{F}$ is given by 
$$X=(iz_1,i\theta_2z_2,\cdots,i\theta_m z_m).$$
For an integer $l\in\{1,\cdots,m-1\}$ and $p\in\{1,\cdots,m-2\}$, we define the vector fields $Y_l$ and $W_p$ on the tangent space of  $\mathbb{S}^{2m-1}$ by the following \cite{HT}
$$Y_l=(0,\cdots,0,-(\sum_{k=l+1}^m|z_k|^2)z_l,|z_l|^2 z_{l+1},\cdots,|z_l|^2z_m),$$ 
and, 
$$W_p=(0,\cdots,0,-(\sum_{k=p+1}^m\theta_k^2|z_k|^2)iz_p, \theta_p\theta_{p+1}|z_p|^2iz_{p+1},\cdots, \theta_p\theta_m|z_p|^2iz_m).$$ 
We also denote by $W_{m-1}$ the vector field on $T\mathbb{S}^{2m-1}$ by 
$$W_{m-1}=(0,\cdots,0,-\theta_m|z_m|^2iz_{m-1},\theta_{m-1}|z_{m-1}|^2 iz_{m}).$$  
It is easy to see that the set $\{X,Y_l, W_p, W_{m-1}\}$ is an orthogonal frame of the tangent space of the sphere for any $l$ and $p.$ Recall now that given an orthonormal frame $\{X/|X|, e_i=Z_i/|Z_i|\}$ of the tangent space of the round sphere for $i=1,\cdots,2m-2$, the norm of O'Neill tensor can be computed as follows
$$|A|^2=\sum_{i,j}|A_{e_i}e_j|^2=\frac{1}{4}\sum_{i,j}|\pi^\perp([e_i,e_j])|^2=\frac{1}{2|X|^2}\sum_{i<j}\frac{1}{|Z_i|^2 |Z_j|^2}|([Z_i,Z_j],X)|^2,$$ 
where $\pi^\perp:TM\rightarrow \R X$ is the projection. On the one hand, a straightforward computation of the norms yields to 
$$|X|^2=|z_1|^2+\sum_{k=2}^m\theta_k^2|z_k|^2.$$ 
Moreover, for any $l$ and $p$, we have
$$|Y_l|^2=|z_l|^2\left(\sum_{k=l+1}^m|z_k|^2\right)\left(\sum_{k=l}^m|z_k|^2\right),$$ 
$$|W_p|^2=|z_p|^2\left(\sum_{t=p+1}^m\theta_t^2|z_t|^2\right)\left(\sum_{s=p}^m\theta_s^2|z_s|^2\right).$$ 
Also we find that, 
$$|W_{m-1}|^2=(\theta_m^2|z_m|^2+\theta_{m-1}^2|z_{m-1}|^2)|z_m|^2|z_{m-1}|^2.$$  
On the other hand, the computation of the Lie brackets yields for any $l$ to 
$$([Y_l,W_l],X)=-2\theta_l|z_l|^2\left(\sum_{s=l+1}^m \theta_s^2|z_s|^2\right)\left(\sum_{k=l}^m|z_k|^2\right),$$ 
and for $l>p,$  
$$([Y_l,W_p],X)=2|z_l|^2\theta_p|z_p|^2\sum_{k=l+1}^m(\theta_l^2-\theta_k^2)|z_k|^2.$$ 
Also, we have that 
$$([Y_{m-1},W_{m-1}],X)=-2|z_{m-1}|^2|z_m|^2\theta_{m-1}\theta_m(|z_{m-1}|^2+|z_m|^2).$$
The other Lie brackets are all equal to zero. Thus the O'Neill tensor is equal to  
\begin{eqnarray*} 
|A|^2&=&\frac{2}{|X|^2}\{\frac{\theta_{m-1}^2\theta_m^2(|z_{m-1}|^2+|z_m|^2)}{\theta_{m-1}^2|z_{m-1}|^2+\theta_m^2|z_m|^2)}\\
&&+\sum_{j=1}^{m-2}\frac{\theta_j^2(\sum_{s=j+1}^m\theta_s^2|z_s|^2)(\sum_{k=j}^m|z_k|^2)}{(\sum_{s=j}^m\theta_s^2|z_s|^2)(\sum_{k=j+1}^m|z_k|^2)}\\
&&+\sum_{j=1}^{m-2}\sum_{i=j+1}^{m-1}\frac{|z_i|^2|z_j|^2(\sum_{k=i+1}^m(\theta_i^2-\theta_k^2)|z_k|^2)^2}{(\sum_{t=j+1}^m\theta_t^2|z_t|^2)(\sum_{s=j}^m\theta_s^2|z_s|^2)(\sum_{k=i+1}^m|z_k|^2)(\sum_{k=i}^m|z_k|^2)}\}.
\end{eqnarray*}
We will now prove that the norm is constant if only if all $\theta_i's$ are equal to $1.$ Indeed, if we evaluate this norm when it corresponds to the cases where $|z_m|\rightarrow 1,\,\, |z_i|\rightarrow 0, i\neq m$ and $|z_{m-1}|\rightarrow 1,\,\, |z_i|\rightarrow 0, i\neq m-1$, we find after identifying that $\theta_{m-1}=\theta_m=\theta.$ The value of the O'Neill tensor corresponding to the case $|z_1|^2=|z_m|^2\rightarrow \frac{1}{2},\,\, |z_i|\rightarrow 0, 2\leq i\leq m-1$ gives that $\theta=1.$ The same computation can be done successively to prove that $\theta_i's$ are equal to one for $i\neq m,m-1$ when considering the case $|z_l|^2=|z_m|^2\rightarrow \frac{1}{2},\,\, |z_i|\rightarrow 0.$ Comparing the lower bound of the inequality in Theorem \ref{thm:1} with the norm of the O'Neill tensor which is equal to $2(m-1),$ we find that the optimality is realized for $\mathbb{S}^5.$\\     

\noindent Next, we will get another pinching condition which doesn't require the positivity of the sectional curvature. We have
\begin{ethm} \label{thm:2}
Under the same condition as in Theorem \ref{thm:1}, we have
$$(q-2)|A|^2\geq {\rm Scal}^M-K^M_1(n-q)(n+q-1)-(p(p-1)+(q-p)(q-p-1))\rho^M_1.$$
\end{ethm}

\noindent {\bf Proof.} The proof is a direct consequence from the fact that
$$\sum_{1\leq i,l\leq q}R^M_{lili}\geq {\rm Scal}_M-K_1^M(n-q)(n+q-1).$$
\hfill$\square$

\noindent The inequality in Theorem \ref{thm:2} is of interest if $$K^M_1(n-q)(n+q-1)+(p(p-1)+(q-p)(q-p-1))\rho^M_1\leq {\rm Scal}^M$$ which with the use of ${\rm Scal}^M\leq K^M_1n(n-1)$ gives that $K^M_1>0.$ 

\section{Foliations with basic harmonic forms}
\setcounter{equation}{0}
In this section, we study the case of Riemannian foliation carrying a basic harmonic form. That is a $p$-form $\alpha$ such that $\Delta_b\alpha=0.$  
\begin{prop} Let $(M,g,\mathcal{F})$ be a Riemannian  manifold carrying a Riemannian foliation. Then, we have
\begin{eqnarray}
2\langle R(\alpha),\alpha\rangle &\geq& -\frac{p-7}{3}\sum_{1\leq i,j\leq q}{\rm Ric}^\nabla_{ij}\langle e_i\lrcorner\alpha,e_j\lrcorner\alpha\rangle+\frac{p-1}{3}\sum_{1\leq i,j,l\leq q} R^M_{lilj}\langle e_i\lrcorner\alpha,e_j\lrcorner\alpha\rangle\nonumber\\
&&-\sum_{1\leq i,j,k,l\leq q} R^M_{ijkl}\langle e_j\wedge e_i\lrcorner\alpha,e_l\wedge e_k\lrcorner\alpha\rangle-\sum_{s=1}^{n-q}\left\{\sum_{i=1}^q|A_{e_i}V_s\lrcorner\alpha|^2+2|\sum_{i=1}^q A_{e_i}V_s\wedge e_i\lrcorner \alpha|^2\right\},\nonumber\\
\label{eq:har}
\end{eqnarray}
for any basic $p$-form $\alpha.$
\end{prop}
{\bf Proof.} For $p=1,$ the inequality is clearly satisfied by \eqref{eq:2}. In order to prove the inequality for $p\geq 2$, we introduce as in \cite{Gro} the operator
$$\mathcal{B}^-\alpha=\frac{1}{(p-2)!}\sum_{i,i_1,\cdots,i_{p-2}}((e_i\wedge e_{i_1}\wedge\cdots \wedge e_{i_{p-2}})\lrcorner\alpha\wedge A_{e_i})\otimes e_{i_1}^*\wedge\cdots\wedge e_{i_{p-2}}^*.$$
The norm of the tensor $\mathcal{B}^-\alpha$ is being defined as the sum
$$|\mathcal{B}^-\alpha|^2=\frac{1}{(p-2)!}\sum_{k,l,i_1,\cdots,i_{p-2}}| (\mathcal{B}^-\alpha)_{kli_{1}\cdots i_{p-2}}|^2.$$
Therefore, we compute
\begin{eqnarray*}
\frac{(p-2)!}{2}|\mathcal{B}^-\alpha|^2&=&\frac{1}{2}\sum_{k,l,i_1,\cdots,i_{p-2}}| (\mathcal{B}^-\alpha)_{kli_{1}\cdots i_{p-2}}|^2\\
&=&\frac{1}{2}\sum_{\substack{i_1,\cdots,i_{p-2}\\i,j,k,l}}\langle((e_i\wedge e_{i_1}\wedge\cdots \wedge e_{i_{p-2}})\lrcorner\alpha\wedge A_{e_i})_{kl}, ((e_j\wedge e_{i_1}\wedge\cdots \wedge e_{i_{p-2}})\lrcorner\alpha\wedge A_{e_j})_{kl}\rangle\\
&=&\frac{1}{2}\sum_{\substack{i_1,\cdots,i_{p-2}\\i,j,k,l}}\langle\alpha_{ii_1\cdots i_{p-2}k}A_{e_i}e_l-\alpha_{ii_1\cdots i_{p-2}l}A_{e_i}e_k,\alpha_{ji_1\cdots i_{p-2}k}A_{e_j}e_l-\alpha_{ji_1\cdots i_{p-2}l}A_{e_j}e_k\rangle\\
&=&\sum_{\substack{i_1,\cdots,i_{p-2}\\i,j,k,l}}\alpha_{iki_1\cdots i_{p-2}}\alpha_{jki_1\cdots i_{p-2}}g(A_{e_i}e_l,A_{e_j}e_l)-\alpha_{iki_1\cdots i_{p-2}}\alpha_{jli_1\cdots i_{p-2}}g(A_{e_i}e_l,A_{e_j}e_k)\\
&=&(p-1)!\sum_{1\leq i,j,l\leq q}\langle e_i\lrcorner\alpha,e_j\lrcorner\alpha\rangle g(A_{e_i}e_l,A_{e_j}e_l)\\
&&-(p-2)!\sum_{1\leq i,j,k,l\leq q}\langle e_k\wedge e_i\lrcorner\alpha,e_l\wedge e_j\lrcorner\alpha\rangle g(A_{e_i}e_l,A_{e_j}e_k).
\end{eqnarray*}
Thus we deduce that
\begin{eqnarray*}
\frac{1}{2}|\mathcal{B}^-\alpha|^2&=&(p-1)\sum_{1\leq i,j,l\leq q}\langle e_i\lrcorner\alpha,e_j\lrcorner\alpha\rangle g(A_{e_i}e_l,A_{e_j}e_l)\\
&&-\sum_{1\leq i,j,k,l\leq q}\langle e_j\wedge e_i\lrcorner\alpha,e_l\wedge e_k\lrcorner\alpha\rangle g(A_{e_i}e_l,A_{e_k}e_j).
\end{eqnarray*}
Plugging now Equations \eqref{eq:21} and \eqref{eq:22} into the above one, we find that
\begin{eqnarray*}
\frac{1}{2}|\mathcal{B}^-\alpha|^2&=&\frac{p-1}{3}\left\{\sum_{1\leq i,j\leq q}{\rm Ric}^\nabla_{ij}\langle e_i\lrcorner\alpha,e_j\lrcorner\alpha\rangle-\sum_{1\leq i,j,l\leq q}R^M_{lilj}\langle e_i\lrcorner\alpha,e_j\lrcorner\alpha\rangle\right\}\\
&&+\sum_{1\leq i,j,k,l\leq q}\left\{-R^\nabla_{ijkl}+R^M_{ijkl}+2g(A_{e_i}e_j,A_{e_k}e_l)-g(A_{e_k}e_i,A_{e_j}e_l)\right\}\langle e_j\wedge e_i\lrcorner\alpha,e_l\wedge e_k\lrcorner\alpha\rangle.
\end{eqnarray*}
Using Equations \eqref{eq:2} and \eqref{eq:4}, we get
\begin{eqnarray*}
\frac{1}{2}|\mathcal{B}^-\alpha|^2&=&\frac{p-1}{3}\left\{\sum_{1\leq i,j\leq q}{\rm Ric}^\nabla_{ij}\langle e_i\lrcorner\alpha,e_j\lrcorner\alpha\rangle-\sum_{1\leq i,j,l\leq q}R^M_{lilj}\langle e_i\lrcorner\alpha,e_j\lrcorner\alpha\rangle\right\}\\
&&+2\langle R(\alpha),\alpha\rangle-2\sum_{1\leq i,j\leq q}{\rm Ric}^\nabla_{ij}\langle e_i\lrcorner\alpha,e_j\lrcorner\alpha\rangle+\sum_{1\leq i,j,k,l\leq q}R^M_{ijkl}\langle e_j\wedge e_i\lrcorner\alpha,e_l\wedge e_k\lrcorner\alpha\rangle\\
&&+2\sum_{s=1}^{n-q}|(\sum_{i=1}^q A_{e_i}V_s\wedge e_i)\lrcorner\alpha|^2-\sum_{1\leq i,j,k,l\leq q}\langle e_j\wedge e_i\lrcorner\alpha,e_l\wedge e_k\lrcorner\alpha\rangle g(A_{e_i}e_l,A_{e_j}e_k).
\end{eqnarray*}
Finally, after the use of Equation \eqref{eq:6}, we find with the help of \eqref{eq:5} that
\begin{eqnarray*}
\frac{1}{2}|\mathcal{B}^-\alpha|^2+|\mathcal{B}^+\alpha|^2&=&\frac{p-7}{3}\sum_{1\leq i,j\leq q}{\rm Ric}^\nabla_{ij}\langle e_i\lrcorner\alpha,e_j\lrcorner\alpha\rangle-\frac{p-1}{3}\sum_{1\leq i,j,l\leq q}R^M_{lilj}\langle e_i\lrcorner\alpha,e_j\lrcorner\alpha\rangle\\&&+2\langle R(\alpha),\alpha\rangle+\sum_{1\leq i,j,k,l\leq q}R^M_{ijkl}\langle e_j\wedge e_i\lrcorner\alpha,e_l\wedge e_k\lrcorner\alpha\rangle\\
&&+\sum_{s=1}^{n-q}\left\{\sum_{i=1}^q|A_{e_i}V_s\lrcorner\alpha|^2+2|(\sum_{i=1}^q A_{e_i}V_s\wedge e_i)\lrcorner\alpha|^2\right\}.
\end{eqnarray*}
Since the l.h.s. of the equality above is non-negative, we finish the proof of the proposition.  
\hfill$\square$\\\\

\noindent We investigate now the case where the form $\alpha$ is a basic harmonic form. We have
\begin{ethm} \label{thm:harmo}
Let $(M,g,\mathcal{F})$ be a compact Riemannian manifold endowed with a Riemannian foliation of codimension $q$. Assume that the normal bundle carries a basic harmonic $p$-form, there exists at least a point $x\in M$ such that
$$(2q+1)|A|^2(x)\geq -\frac{p-7}{3}{\rm Scal}^\nabla(x)+(\frac{p-1}{3})q(q-1)K^M_0(x)-2(p(p-1)+(q-p)(q-p-1))\rho_1^M(x),$$
where $2\leq p\leq q-2.$
\end{ethm}
{\bf Proof.} As in the proof of Theorem \ref{thm:1}, we use Inequality \eqref{eq:9} in order to deduce that 
\begin{eqnarray*}
2\langle R(\alpha),\alpha\rangle&\geq& -\frac{p-7}{3}\sum_{1\leq i,j\leq q}{\rm Ric}^\nabla_{ij}\langle e_i\lrcorner\alpha,e_j\lrcorner\alpha\rangle+\frac{p-1}{3}\sum_{1\leq i,j,l\leq q} R^M_{lilj}\langle e_i\lrcorner\alpha,e_j\lrcorner\alpha\rangle\\&&-\sum_{1\leq i,j,k,l\leq q} R^M_{ijkl}\langle e_j\wedge e_i\lrcorner\alpha,e_l\wedge e_k\lrcorner\alpha\rangle-(2q+1)\sum_{\substack{1\leq i\leq q\\1\leq s\leq n-q}}|A_{e_i}V_s\lrcorner\alpha|^2
\end{eqnarray*}
for any basic $p$-form $\alpha$. Applying the above inequality for the $(q-p)$-form $*_b\alpha$ and then summing the two equations, we get by using \eqref{eq:cour} and \eqref{eq:8}
\begin{eqnarray*}
2(\langle R(\alpha),\alpha\rangle+\langle R(*_b\alpha),*_b\alpha\rangle)&\geq& -\frac{p-7}{3}{\rm Scal}^\nabla |\alpha|^2+\frac{p-1}{3}q(q-1)K^M_0|\alpha|^2\\
&&-2(p(p-1)+(q-p)(q-p-1))\rho_1^M|\alpha|^2\\&&-(2q+1)|A|^2|\alpha|^2.
\end{eqnarray*}
If the form $\alpha$ is now harmonic, i.e. $d_b\alpha=\delta_b\alpha=0$, then the twisted derivative is equal to $\tilde{d}_b\alpha=-\frac{1}{2}\kappa\wedge\alpha$ and its adjoint is $\tilde{\delta}_b\alpha=-\frac{1}{2}\kappa\lrcorner\alpha.$ Thus
$$\int_M\langle \tilde{\Delta}_b\alpha,\alpha\rangle v_g=\int_M (|\tilde{d}_b\alpha|^2+|\tilde{\delta}_b\alpha|^2)v_g=\frac{1}{4}\int_M|\kappa|^2|\alpha|^2v_g.$$
This implies by the transverse Bochner-Weitzenb\"ock formula,
that $\int_M\langle R(\alpha),\alpha\rangle v_g\leq 0.$ Since the basic Hodge operator commutes with
the twisted Laplacian \cite{HR}, the same inequality
holds for $*_b\alpha$. Thus, we get the required
inequality.
\hfill$\square$\\\\


\end{document}